\documentclass[10pt,draft,reqno]{amsart}
    \makeatletter
    \def\section{\@startsection{section}{1}%
    \z@{.7\linespacing\@plus\linespacing}{.5\linespacing}%
    {\bfseries
    \centering
    }}
    \def\@secnumfont{\bfseries}
    \makeatother
\newtheorem{theorem}{Theorem}[section]
\newtheorem{lemma}[theorem]{Lemma}
\newtheorem{proposition}[theorem]{Proposition}
\newtheorem{corollary}[theorem]{Corollary}
\theoremstyle{definition}
\newtheorem{definition}[theorem]{Definition}
\newtheorem{example}[theorem]{Example}
\theoremstyle{remark}
\newtheorem{remark}[theorem]{Remark}
\numberwithin{equation}{section} 

\newcommand{\ombar}{\overline{\Omega}}
\newcommand{\po}{\partial\Omega}
\newcommand{\me}{\mathcal E}
\newcommand{\mf}{\mathcal F}
\newcommand{\mfm}{\mathcal F^{\mu}}
\newcommand{\mem}{\mathcal E^{\mu}}
\newcommand{\Huntild}{\widetilde{H}^1(\Omega)}
\newcommand{\pmu}{\mathcal P_t^{\mu}}
\newcommand{\ma}{\mathcal A}

\begin{document}

\title[]
{Probabilistic Solution of the General Robin Boundary Value Problem  on Arbitrary Domains}

\author{ AKHLIL KHALID }

\address{Khalid Akhlil: Department of Mathematics, Ibno Zohr university, Agadir, Morocco.}

\email{khalid.akhlil@uni-ulm.de}

\thanks{}

\thanks{}



\date{}

\dedicatory{}

\commby{}


\begin{abstract}

Using a capacity approach, and the theory of measure's
perturbation of Dirichlet forms, we give the probabilistic
representation of the General Robin boundary value problems on an
arbitrary domain $\Omega$, involving smooth measures, which give
arise to a new process obtained by killing the  general
reflecting Brownian motion at a random time. We obtain some
properties of the semigroup directly from its probabilistic
representation, and some convergence theorems, and also a probabilistic interpretation
of the phenomena occurring on the boundary.
\end{abstract}

\maketitle

\section{Introduction}

The classical Robin boundary conditions on a smooth domain
$\Omega$ of $\mathbb R^N$ ($N\geq 0$), is giving by :
\begin{equation}\frac{\partial
u}{\partial\nu}+\beta u=0\quad\text{on }\po,
\end{equation}
where $\nu$ is the outward normal vector field on the boundary
$\po$, and $\beta$ a positive bounded Borel measurable function
defined on $\po$.

The probabilistic treatment of Robin boundary value problems has
been considered by many authors \cite{SU, P, R,Re}. The first two
authors considered bounded $C^3-$domains since the third
considered bounded domains with Lipschitz boundary, and \cite{Re}
was concerned with $C^3-$ domains but with smooth measures instead
of $\beta$. If one want to generalize the probabilistic treatment
to a general domains, a difficulty arise when we try to get a
diffusion process representing Neumann boundary conditions.

In fact, the Robin boundary conditions (1.1) are nothing but a
perturbation of $\frac{\partial}{\partial\nu}$, which represent
Neumann boundary conditions, by the measure $\mu=\beta.\sigma$,
where $\sigma$ is the surface measure. Consequently, the
associated diffusion process is the reflecting Brownian motion
killed by a certain additive functional, and the semigroup
generated by the Laplacian with classical Robin boundary
conditions is then giving by:
\begin{equation}\pmu f(x)=E_x[f(X_t)e^{-\int_0^t\beta(X_s)dL_s}]
\end{equation}
where $(X_t)_{t\geq 0}$ is a reflecting Brownian motion (RBM), and
$L_t$ is the boundary local time, which corresponds to $\sigma$ by
Revuz correspondence. It is clear that the smoothness of the
domain $\Omega$ in classical Robin boundary value problem, follows
the smoothness of the domains where RBM is constructed( see
\cite{BBC2, BH2, BH, CFW, F, H3} and references therein for more
details about RBM).

In \cite{BH2}, the RBM is defined to be the Hunt process
associated with the form $(\me,\mf)$ defined on $L^2(\Omega)$ by:
\begin{equation*}
\me(u,v)=\int_{\Omega}\nabla u\nabla v dx\quad,\forall
u,v\in\mf=H^1(\Omega)
\end{equation*}
where $\Omega$ is assumed to be bounded with Lipschitz boundary so
that the Dirichlet form $(\me,\mf)$ is regular. If $\Omega$ is an
arbitrary domain, then the Dirichlet form need not to be Regular,
and to not loose the generality we consider $\mf=\Huntild$, the
closure of $H^1(\Omega)\cap C_c(\overline{\Omega})$ in
$H^1(\Omega)$. The domain $\widetilde{H}^1(\Omega)$ is so defined
to insure the Dirichlet form $(\me,\mf)$ to be regular.

Now, if we perturb the Neumann boundary conditions by Borel
positive measure \cite{AW1, AW2, Wa}, we get the Dirichlet form
$(\mem,\mfm)$ defined on $L^2(\Omega)$ by:
\begin{equation}
\mem(u,v)=\int_{\Omega}\nabla u\nabla v
dx+\int_{\po}\widetilde{u}\widetilde{v}d\mu\quad,\forall
u,v\in\mfm=\Huntild\cap L^2(\po,d\mu)
\end{equation}
In the case of $\mu=\beta.\sigma$($\Omega$ bounded with Lipschitz
boundary), (1.3) is the form associated with Laplacian with classical
Robin boundary conditions and (1.2) gives the associated
semigroup. In the case of an arbitrary domain $\Omega$ we make use
of the theory of measure's perturbation of Dirichlet forms, see
e.g. \cite{AM1, AM2, BM, FOT, K, Si, S, SV, V1, V2}.

More specifically, we adapt the potential theory, and associated
stochastic analysis to our context, this is the subject of section
2. In section 3, we focus on the diffusion process $(X_t)_{t\geq
0}$ associated with the regular Dirichlet form $(\me,\Huntild)$.
We apply a decomposition theorem of additive functionals to write
$X_t$ in the form $X_t=x+B_t+N_t$, we prove that the additive
functional $N_t$ is supported by $\po$, and we investigate when it
is of bounded variations.

In section 3 we get the probabilistic representation of the
semigroup associated with (1.3), and we prove that it is
sandwiched between the semigroup generated by the Laplacian with
Dirichlet boundary conditions, and that of Neumann ones. In
addition, we prove some convergence theorems,
and we give a probabilistic interpretation of the
phenomena occurring on the boundary.

\section{Preliminaries and notations}

The aim of this section is to adapt the potential theory, and the
stochastic analysis for application to our problem. More
precisely, it concerns the notion of relative capacity, smooth
measures, and its corresponding additive functionals. This section
relies heavily on the book of Fukushima\cite{FOT}, particulary
chapter 2 and 5, and the paper \cite{AW1}. Throughout \cite{FOT},
the form $(\me,\mf)$ is a regular Dirichlet form on $L^2(X,m)$,
where $X$ is a locally compact separable metric space, and $m$ a
positive Radon measure on $X$ with $supp[m]=X$.

For our purposes we take $X=\ombar$, where $\Omega$ is an Euclidean domain of $\mathbb R^N$, and the measure $m$ on the
$\sigma-$algebra $\mathcal B(X)$ is given by
$m(A)=\lambda(A\cap\Omega)$ for all $A\in\mathcal B(X)$ with
$\lambda$ the Lebesgue measure, it follows that
$L^2(\Omega)=L^2(X,\mathcal B(X),m)$, and we define a regular
Dirichlet form $(\me,\mf)$ on $L^2(\Omega)$ by:

$$\me(u,v)=\int_{\Omega}\nabla u\nabla v dx\quad,
\mf=\widetilde{H^1}(\Omega)$$where
$\widetilde{H}^1(\Omega)=\overline{H^1(\Omega)\cap
C_c(\overline{\Omega})}^{H^1(\Omega)}$. The domain
$\widetilde{H}^1(\Omega)$ is so defined to insure the Dirichlet
form $(\me,\mf)$ to be regular, instead of $\mf= H^1(\Omega)$
which make the form not regular in general, but if $\Omega$ is
bounded open set with Lipschiz boundary, then
$\widetilde{H}^1(\Omega)=H^1(\Omega)$.

We denote for any $\alpha>0:\text{
}\me_{\alpha}(u,v)=\me(u,v)+\alpha(u,v)_m,\text{ }\forall
u,v\in\mf$.

\subsection{ Relative Capacity} The relative capacity is introduced
in a first time in \cite{AW1} to study the Laplacian with general
Robin boundary conditions on arbitrary domains. It is a special
case of the capacity associated with a regular Dirichlet form as
described in chapter 2 of \cite{FOT}. It seems to be an efficient
tool to analyse the phenomena occurring on the boundary $\po$ of
$\Omega$.

The relative capacity which we denote by
$Cap_{\overline{\Omega}}$ is defined on a subsets of
$\overline{\Omega}$ by: For $A\subset\overline{\Omega}$ relatively
open (i.e. open with respect to the topology of
$\overline{\Omega}$) we set:$$Cap_{\overline{\Omega}}(A):=\text{inf}\{\me_1(u,u):u\in\widetilde{H}^1(\Omega):u\geq 1\text{ a.e on
}A\}$$

And for arbitrary $A\subset\overline{\Omega}$, we set:
$$Cap_{\overline{\Omega}}(A):=\text{inf}\{Cap_{\overline{\Omega}}(B):B
\text{ relatively open }A\subset B\subset\overline{\Omega}\}$$

A set $N\subset\overline{\Omega}$ is called a relatively polar if
$Cap_{\overline{\Omega}}(N)=0$.

The relative capacity (just as a cap) has the properties of a
capacity as described in \cite{FOT}. In particular, $cap_{\ombar}$
is also an outer measure (But not a Borel measure) and a Choquet
Capacity.

A statement depending on $x\in A\subset\overline{\Omega}$ is said
to hold relatively quasi-everywhere (r.q.e.) on $A$, if there
exist a relatively polar set $N\subset A$ such that the statement
is true for every $x\in A\setminus N$.

Now we may consider functions in $\widetilde{H}^1(\Omega)$ as
defined on $\overline{\Omega}$, and we call a function
$u:\overline{\Omega}\rightarrow\mathbb R$ relatively
quasi-continuous (r.q.c.) if for every $\epsilon>0$ there exists a
relatively open set $G\subset\overline{\Omega}$ such that
$Cap_{\overline{\Omega}}(G)<\epsilon$ and
$u|_{\overline{\Omega}\setminus G}$ is continuous.

It follows \cite{Wa}, that for each $u\in\Huntild$ there exists a
relatively quasi-continuous function
$\widetilde{u}:\ombar\rightarrow\mathbb R$ such that
$\widetilde{u}(x)=u(x)$ $m-$a.e. This function is unique
relatively quasi-everywhere. We call $\widetilde{u}$ the
relatively quasi-continuous representative of $u$.

For more details, we refer the reader to \cite{AW1, Wa}, where the
relative capacity is investigated, as well as its relation to the
classical one. A description of the space $H^1_0(\Omega)$ is term
of relative capacity is also given, namely:
\begin{equation}H^1_0(\Omega)=\{u\in\Huntild:\widetilde{u}(x)=0\text{ r.q.e. on }\po\}
\end{equation}

\subsection{Smooth measures}

All families of measures on $\po$ defined in this subsection, was
originally defined on $X$ \cite{FOT}, and then in our settings on
$X=\ombar$, as a special case. We reproduce the same definitions,
and most of their properties on $\po$, as we deal with measures
concentrated on the boundary of $\Omega$ for our approach to Robin
boundary conditions involving measures. There is three families of
measures, as we will see in the sequel, the family $S_0$, $S_{00}$
and $S$. We put $\po$ between brackets to recall our context, and
we keep in mind that the same things are valid if we put $\Omega$
or $\ombar$ instead of $\po$.

Let $\Omega\subset\mathbb R^N$ be open. A positive Radon measure
$\mu$ on $\partial\Omega$ is said to be of finite energy integral
if $$\int_{\partial\Omega}|v(x)|\mu(dx)\leq C\sqrt{\mathcal
E_1(v,v)}\quad,v\in\mathcal F\cap C_c(\overline{\Omega})$$for some
positive constant $C$. A positive Radon measure on
$\partial\Omega$ is of finite energy integral if and only if there
exists, for each $\alpha>0$, a unique function
$U_{\alpha}\mu\in\mathcal F$ such that $$\mathcal
E_{\alpha}(U_{\alpha}\mu,v)=\int_{\partial\Omega}v(x)\mu(dx)$$ We
call $U_{\alpha}\mu$ an $\alpha-$potential.

We denote by $S_0(\partial\Omega)$, the family of all positive
Radon measures of finite energy integral.

\begin{lemma} Each measure in $S_0(\partial\Omega)$ charges no set of zero
relative capacity.\end{lemma}

Let us consider a subset $S_{00}(\po)$ of $S_0$ defined
by:$$S_{00}(\po)=\{\mu\in
S_0(\po):\mu(\po)<\infty,||U_1\mu||_{\infty}<\infty\}$$

\begin{lemma}
For any $\mu\in S_0(\po)$, there exist an increasing sequence
$(F_n)_{n\geq 0}$ of compact sets of $\po$ such that:
$$1_{F_n}.\mu\in S_{00}(\po)\quad,n=1,2,...$$
$$Cap_{\ombar}(K\setminus F_n)\longrightarrow 0 ,n\rightarrow
+\infty \text{ for any compact set }K\subset\po$$

\end{lemma}

We note that $\mu\in S_0(\po)$ vanishes on
$\po\setminus\cup_{n}F_n$ for the sets $F_n$ of the Lemma 2.2,
because of the Lemma 2.1 .

We now turn to a class of measures $S(\po)$ larger than
$S_0(\po)$. Let us call a (positive) Borel measure $\mu$ on $\po$
smooth if it satisfies the following conditions:

- $\mu$ charges no set of zero relative capacity.

- There exist an increasing sequence $(F_n)_{n\geq 0}$ of closed
sets of $\po$ such that: \begin{equation}
\mu(F_n)<\infty\quad,n=1,2,...\end{equation}

\begin{equation} \lim_{n\rightarrow +\infty}Cap_{\ombar}(K\setminus F_n)=0\text{
for any compact } K\subset\po
\end{equation}
Let us note that $\mu$ then satisfies
\begin{equation}\mu(\po\setminus\cup_n F_n)=0\end{equation}

An increasing sequence $(F_n)$ of closed sets satisfying condition
(2.3) will be called a generalized nest, if further each $F_n$ is
compact, we call it a generalized compact nest.

We denote by $S(\po)$ the family of all smooth measures. The class
$S(\po)$ is quiet large and it contains all positive Radon measure
on $\po$ charging no set of zero relative capacity. There exist
also, by Theorem 1.1 \cite{AM2} a smooth measure $\mu$ on $\po$ (
hence singular with respect to $m$) "nowhere Radon" in the sense
that $\mu(G)=\infty$ for all non-empty relatively open subset $G$
of $\po$ (See Example 1.6\cite{AM2}).

The following Theorem, say that, any measure in $S(\po)$ can be
approximated by measures in $S_0(\po)$ and in $S_{00}(\po)$ as
well.

\begin{theorem} The following conditions are equivalent for a positive Borel
measure $\mu$ on $\po$:

(i) $\mu\in S(\po)$.

(ii) There exists a generalized nest $(F_n)$ satisfying $(2.4)$
and $1_{F_n}.\mu\in S_0(\po)$ for each $n$.

(iii) There exists a generalized compact nest $(F_n)$ satisfying
$(2.4)$ and $1_{F_n}.\mu\in S_{00}(\po)$ for each $n$.
\end{theorem}

\subsection{Additive functionals}

Now we turn our attention to the correspondence between smooth
measures and additive functionals, known as Revuz correspondence.
As the support of an additive functional is the quasi-support of
its Revuz measure, we restrict our attention, as for smooth
measures, to additive functionals supported by $\po$. Recall that
as the Dirichlet form $(\me,\mf)$ is regular, then there exists a
Hunt process $M=(\Xi,X_t,\xi,P_x)$ on $\ombar$ which is
$m-$symmetric and associated with it.

\begin{definition}A function
$A:[0,+\infty[\times\Xi\rightarrow[-\infty,+\infty]$ is said to be
an Additive functional (AF) if:

1) $A_t$ is $\mathcal F_t-$measurable.

2) There exist a defining set $\Lambda\in\mathcal F_{\infty}$ and
an exceptional set $N\subset \po$ with $cap_{\ombar}(N)=0$ such
that $P_x(\Lambda)=1$, $\forall x\in \po\setminus N$,
$\theta_t\Lambda\subset\Lambda$, $\forall t>0$;
$\forall\omega\in\Lambda, A_0(\omega)=0$; $|A_t(\omega)|<\infty$
for $t<\xi$. $A_{.}(\omega)$ is right continuous and has left
limit, and $ A_{t+s}(\omega)=A_t(\omega)+A_s(\theta_t\omega)\text{
}s,t\geq 0$\end{definition}

An additive functional is called positive continuous (PCAF) if, in
addition, $A_t(\omega)$ is nonnegative and continuous for each
$\omega\in\Lambda$. The set of all PCAF's on $\po$ is denoted
$\ma_c^{+}(\po)$.

Two additive functionals $A^1$ and $A^2$ are said to be equivalent
if for each $t>0$, $P_x(A^1_t=A^2_t)=1\text{ r.q.e }x\in \ombar$.

We say that $A\in\ma_c^+(\po)$ and $\mu\in S(\po)$ are in the
Revuz correspondence, if they satisfy, for all $\gamma-$excessive
function $h$, and $f\in\mathcal B_+(\ombar)$,  the
relation:$$\lim_{t\searrow
0}\frac{1}{t}E_{h.m}\left[\int_0^tf(X_s)dA_s\right]=\int_{\po}h(x)(f.\mu)(dx)$$
The family of all equivalence classes of $\ma_c^+(\po)$ and the
family $S(\po)$ are in one to one correspondence under the Revuz
correspondence. In this case, $\mu\in S(\po)$ is called the Revuz
measure of $A$.

\begin{example}
We suppose $\Omega$ to be bounded with Lipschitz boundary.
We have \cite{P}:
$$\lim_{t\searrow
0}\frac{1}{t}E_{h.m}\left[\int_0^tf(X_s)dL_s\right]=\frac{1}{2}\int_{\partial\Omega}h(x)f(x)\sigma(dx)$$
where $L_t$ is the boundary local time of the reflecting Brownian
motion on $\ombar$. It follows that $\frac{1}{2}\sigma$ is the
Revuz measure of $L_t$ .
\end{example}

In the following we give some facts useful in the proofs of our
main results. We set:

$$U_A^{\alpha}f(x)=E_x[\int_0^{\infty}e^{-\alpha t}f(X_t)d A_t]$$
$$R^A_{\alpha}f(x)=E_x[\int_0^{\infty}e^{-\alpha t}e^{-A_t}f(X_t)d
t]$$
$$R_{\alpha}f(x)=E_x[\int_0^{\infty}e^{-\alpha t}f(X_t)d
t]$$

\begin{proposition}
Let $\mu\in S_0(\po)$ and $A\in\ma_c^+(\po)$ the corresponding
PCAF. For $\alpha>0$, $f\in\mathcal B_b^+$, $U^{\alpha}_A$ is a
relatively quasi-continuous version of $U_{\alpha}(f.\mu)$.
\end{proposition}

\begin{proposition}
Let $A\in\ma_c^+(\po)$, and $f\in\mathcal B_b^+$, then
$R^A_{\alpha}$ is relatively quasi-continuous and
$$R_{\alpha}^Af-R_{\alpha}f+U_A^{\alpha}R_{\alpha}^Af=0$$
\end{proposition}

In general, the support of an AF $A$ is defined by
$$supp[A]=\{x\in X\setminus
N: P_x(R=0)=1\}$$ where $R(\omega)=inf\{t>0:A_t(\omega)\neq 0\}$

\begin{theorem}
The support of $A\in\ma_c^+(\po)$ is the relative quasi-support of
its Revuz measure.

\end{theorem}

In the following we give a well known theorem of decomposition of
additive functionals of finite energy. We will apply it to get a
decomposition of the diffusion process associated with
$(\me,\mf)$.

\begin{theorem}
For any $u\in\mf$, the AF
$A^{[u]}=\widetilde{u}(X_t)-\widetilde{u}(X_0)$ can be expressed
uniquely as
\begin{equation}\widetilde{u}(X_t)-\widetilde{u}(X_0)=M^{[u]}+N^{[u]}\end{equation}
where $M_t^{[u]}$ is a martingale additive functional of finite
energy and $N_t^{[u]}$ is a continuous additive functional of zero
energy.
\end{theorem}

A set $\sigma(u)$ is called the $(0)-$spectrum of $u\in\mf$, if
$\sigma(u)$ is the complement of the largest open set $G$ such
that $\me(u,v)$ vanishes for any $v\in\mf\cap\mathcal C_0(X)$ with
$supp[v]\subset G$. The following Theorem means that :
$supp[N^{[u]}]\subset\sigma(u)$, $\forall u\in\mf$.

\begin{theorem}
For any $u\in\mf$, the CAF $N^{[u]}$ vanishes on the complement of
the spectrum $F=\sigma(u)$ of $u$ in the following
sense:$$P_x(N_t^{[u]}=0:\forall t<\sigma_F)=1\text{ r.q.e }x\in
X$$
\end{theorem}

\section{General Reflecting Brownian Motion}

Now we turn our attention to the process associated with the
regular Dirichlet form $(\me,\mf)$ on $L^2(\Omega)$ defined by:

\begin{equation}\me(u,v)=\int_{\Omega}\nabla u\nabla v dx\quad,
\mf=\widetilde{H^1}(\Omega)\end{equation}

Due to the Theorem of Fukushima(1975), there is a Hunt process
$(X_t)_{t\geq 0}$ associated with it. In addition, $(\me,\mf)$ is
local, thus the Hunt process is in fact a diffusion process (i.e.
A strong Markov process with continuous sample paths). The
diffusion process $M=(X_t,P_x)$ on $\ombar$ is associated with the
the form $\me$ in the sense that the transition semigroup
$p_tf(x)=E_x[f(X_t)]$, $x\in\ombar$ is a version of the
$L^2-$semigroup $\mathcal P_tf$ generated by $\me$ for any
nonnegative $L^2-$function $f$.

$M$ is unique up to set of zero relative capacity.

\begin{definition}
We call the diffusion process on $\ombar$ associated with
$(\me,\mf)$ the \emph{General reflecting Brownian motion}.
\end{definition}

The process $X_t$ is so named to recall the standard reflecting
Brownian motion in the case of bounded smooth $\Omega$, as the
process associated with $(\me,H^1(\Omega))$. Indeed, when $\Omega$
is bounded with Lipschitz boundary we have that
$\Huntild=H^1(\Omega)$, and by \cite{BH2} the reflecting Brownian
motion $X_t$ admits the following Skorohod representation:

\begin{equation}X_t=x+W_t+\frac{1}{2}\int_0^t\nu(X_s)dL_s,\end{equation}where
$W$ is a standard $N-$dimensional Brownian motion, $L$ is the
boundary local(continuous additive functional) associated with
surface measure $\sigma$ on $\po$, and $\nu$ is the inward unit
normal vector field on the boundary.

For a general domains, the form $(\me,H^1(\Omega))$ need not to be
regular. Fukuchima \cite{F} constructed the reflecting brownian
motion on a special compactification of $\Omega$, the so called
Kuramuchi compactification. In \cite{BH2} it is shown that if
$\Omega$ is a bounded Lipschitz domain, then the Kuramochi
compactification of $\Omega$ is the same as Euclidean
Compactification. Thus for such domains, the reflecting Brownian
motion is a continuous process who does live on the set $\ombar$.

Now,we apply a general decomposition theorem of additive
functionals to our process $M$, in the same way as in \cite{BH2}.
According to Theorem 2.9 the continuous additive functional
$\widetilde{u}(X_t)-\widetilde{u}(X_0)$ can be decomposed as
follows:
$$\widetilde{u}(X_t)-\widetilde{u}(X_0)=M_t^{[u]}+N_t^{[u]}$$where
$M_t^{[u]}$ is a martingale additive functional of finite energy
and $N_t^{[u]}$ is a continuous additive functional of zero
energy.

Since $(X_t)_{t\geq 0}$ has continuous sample paths, $M_t^{[u]}$
is a continuous martingale whose quadratic variation process is:

\begin{equation}<M^{[u]},M^{[u]}>_t=\int_0^t|\nabla
u|^2(X_s)ds\end{equation}

Instead of $u$ we take coordinate function $\phi_i(x)=x_i$. We
have $$X_t=X_0+M_t+N_t$$

We claim that $M_t$ is a Brownian motion with respect to the
filtration of $X_t$. To see that, we use L\'evys criterion. This
follows immediately from (3.2), which became in the case of
coordinate function: $$<M^{[\phi_i]},M^{[\phi_i]}>=\delta_{ij}t$$

Now we turn our attention to the additive functional $N_t$. Two
natural questions need to be answered. The first is, where is the
support of $N_t$ located, and the second concern the boundedness
of its total variation.

For the first question we claim the following:

\begin{proposition}
The additive functional $N_t$ is supported by $\po$.
\end{proposition}

\begin{proof}

Following Theorem 2.10, we have that $supp[N_t]\subset \sigma
(\phi)$, where $\sigma(\phi)$ is the $(0)-$spectrum of $\phi$,
which means the complement of the largest open set $G$ such that
$\mathcal E(\phi_i,v)=0$ for all $v\in\mathcal F\cap
C_c(\overline{\Omega})$ with $supp[v]\subset G$.
\bigskip

$\underline{\text{\textit{Step 1}}}$: If $\Omega$ is smooth(
Bounded with Lipschitz boundary, for example), then we have:
$$\mathcal E(\phi_i,v)=-\int_{\po}v.n_i d\sigma$$

Then, $\mathcal E(\phi_i,v)=0$ for all $v\in\mathcal F\cap
C_c(\overline{\Omega})$ with $supp[v]\subset \Omega$. We can then
see that the largest $G$ is $\Omega$. Consequently
$\sigma(\phi)=\ombar\setminus\Omega$, and then
$\sigma(\phi)=\partial\Omega$.

$\underline{\text{\textit{Step 2}}}$: If $\Omega$ is
arbitrary, then we take an increasing sequence of subset of
$\Omega$ such that $\bigcup_{n=0}^{\infty}\Omega_n=\Omega$. Define
the family of Dirichlet forms $(\mathcal E_{\Omega_n},\mathcal
F_{\Omega_n})$ to be the parts of the form $(\me,\mf)$ on each
$\Omega_n$ as defined in section 4.4 of \cite{FOT}. By Theorem
4.4.5 in the same section, we have that $\mathcal
F_{\Omega_n}\subset \mathcal F$ and $\mathcal
E_{\Omega_n}=\mathcal E$ on $\mathcal F_{\Omega_n}\times\mathcal
F_{\Omega_n}$. We have that $\Omega_n$ is the largest open set
such that $\me_{\Omega_n}(\phi_i,v)=0$ for all $v\in\mathcal
F_{\Omega_n}\cap C_c(\overline{\Omega_n})$. By limit, we get the
result.

\end{proof}

The interest of the question of boundedness of total variation of
$N_t$ appears when one need to study the semimartingale property
and the Skorohod equation of the process $X_t$ of type 3.2. Let
$|N|$ be the total variation of $N_t$, i.e.,
$$|N|_t=sup\sum_{i=1}^{n-1}|N_{t_i}-N_{t_{i-1}}|.$$where the
supremum is taken over all finite partition $0=t_0<t_1<...<t_n=t$,
and $|.|$ denote the Euclidian distance. If $|N|$ is bounded, then
we have the following expression:

$$N_t=\int_0^t\nu_s d |N|_s$$where $\nu$ is a process such that
$|\nu|_s=1$ for $|N|-$almost all $s$.

According to \S 5.4. in \cite{FOT}, we have the following result:

\begin{theorem}
Assume that $\Omega$ is bounded, and that the following inequality
is satisfied:

\begin{equation}
\left|\int_{\Omega}\frac{\partial v}{\partial x_i}dx\right|\leq
C||v||_{\infty}\quad,\forall v\in\Huntild\cap C_b(\ombar)
\end{equation}for some constant $C$. Then, $N_t$ is of bounded
variation.
\end{theorem}

A bounded set verifying (2.3) is called \emph{strong Caccioppoli
set}. This notion is introduced in \cite{CFW}, and is a purely
measure theoretic notion. An example of this type of sets are
bounded sets with Lipschitz boundary.

\begin{theorem}
If $\Omega$ is a Caccioppoli set, then there exist a finite signed
smooth measure $\nu$ such that:
\begin{equation}
\int_{\Omega}\frac{\partial v}{\partial
x_i}dx=-\int_{\po}vd\mu\quad ,\forall v\in\Huntild\cap
C_b(\ombar).
\end{equation}
and $\nu=\nu^1-\nu^2$ is associated with the CAF
$-N_t=-A^1_t+A^2_t$ with the Revuz correspondence. Consequently
$\nu$ charges no set of zero relative capacity.
\end{theorem}

To get a Skorohod type representation, we set:
\begin{equation}
\begin{split}
\nu&=\sum_{i=1}^N|\mu_i|\\
\phi_i  & =\frac{d\mu_i}{d\nu}\quad i=1,...,N
\end{split}
\end{equation}
We define the measure $\sigma$ on $\po$ by:
\begin{equation}
\sigma(dx)=2\left(\sum_{i=1}^N|\phi_i(x)|^2\right)^{\frac{1}{2}}\nu(dx)
\end{equation}
and the vector of length 1 at $x\in\po$ by:

$$ n_i(x)= \left\{
       \begin{array}{ll}
           \frac{\phi_i(x)}{\left(\sum_{i=1}^N|\phi_i(x)|^2\right)^{\frac{1}{2}}}  &  if\quad \sum_{i=1}^N|\phi_i(x)|^2>0; \\
           0 &
           if\quad \sum_{i=1}^N|\phi_i(x)|^2=0
         \end{array}
         \right.$$

Thus, $\mu_i(dx)=\frac{1}{2}n_i(x)\sigma(dx)\quad,i=1,..,N$.

Then $$N_t=\int_0^tn(X_s)dL_s$$ where $L$ is the PCAF associated
with $\frac{1}{2}\sigma$.

\begin{theorem}If $\Omega$ is a Caccioppoli set, then for r.q.e
$x\in\ombar$, we have:

$$X_t=x+B_t+\int_0^tn(X_s)dL_s.$$where $B$ is an $N-$dimensional Brownian
motion, and $L$ is a PCAF associated by the Revuz correspondence
to the measure $\frac{1}{2}\sigma$.

\end{theorem}

\begin{remark}
The above theorem can be found in \cite{F} and \cite{F2}. In particular Fukushima proves an equivalence between the property of Caccioppoli sets and the Skorohod representation.
\end{remark}

\section{Probabilistic solution to general Robin boundary value problem}

This section is concerned with the probabilistic representation to
the semigroup generated by the Laplacian with general Robin
boundary conditions, which is, actually, obtained by perturbing
the Neumann boundary conditions by a measure. We start with the
Regular Dirichlet form defined by (3.1), which we call always as
the Dirichlet form associated with the Laplacian with Neumann
boundary conditions.

Let $\mu$ be a positive Radon measure on
$\partial\Omega$ charging no set of zeo relative capacity.
Consider the perturbed Dirichlet form $(\mathcal E_{\mu},\mathcal
F_{\mu})$ on $L^2(\Omega)$ defined by:
$$\mathcal F_{\mu}=\mathcal F \cap L^2(\partial\Omega,\mu)$$
$$\mathcal E_{\mu}(u,v)=\mathcal E(u,v)+\int_{\partial\Omega}uv
d\mu\quad u,v\in\mathcal F_{\mu}$$

We shall see in the following theorem that the transition
function:$$\mathcal P_t^{\mu}f(x)=E_x[f(X_t)e^{-A_t^{\mu}}]$$ is
associated with $(\mathcal E_{\mu},\mathcal F_{\mu})$, where
$A_t^{\mu}$ is a positive additive functional whose Revuz measure
is $\mu$, note that the support of the AF is the same as the
relative quasi-support of its Revuz measure.

\begin{proposition}
$\mathcal P_t^{\mu}$ is a strongly continuous semigroup on
$L^2(\Omega)$.
\end{proposition}
\begin{proof}
The proof of the above Proposition can be found in \cite{AM1}.
\end{proof}

\begin{theorem}Let $\mu$ be a positive Radon measure on $\partial\Omega$
charging no set of zero relative capacity and $(A_t^{\mu})_{t\geq
0}$ be its associated PCAF of $(X_t)_{t\geq 0}$. Then $\mathcal
P_t^{\mu}$ is the strongly continuous semigroup associated with
the Dirichlet form $(\mathcal E_{\mu},\mathcal F_{\mu})$ on
$L^2(\Omega)$.
\end{theorem}

\begin{proof}

To prove that $\mathcal P_t^{\mu}$ is associated with the
Dirichlet form $(\mathcal E_{\mu},\mathcal F_{\mu})$ on
$L^2(\Omega)$ it suffices to prove the assertion
\begin{equation}
\quad R_{\alpha}^A f\in\mathcal F^{\mu}\quad,\mathcal
E_{\alpha}^{\mu}(R_{\alpha}^A,u)=(f,u)\quad,f\in
L^2(\Omega,m),u\in\mathcal F^{\mu}
\end{equation}
Since $\quad ||R_{\alpha}^Af||_{L^2(\Omega)}\leq
||R_{\alpha}f||_{L^2(\Omega)}\leq\frac{1}{\alpha}||f||_{L^2(\Omega)}$,
we need to prove $(4.2)$ only for bounded $f\in L^2(\Omega)$. We
first prove that $(4.2)$ is valid when $\mu\in S_{00}(\po)$.
According to the Proposition 2.7 we have
$$R_{\alpha}^A
f-R_{\alpha}f+U_A^{\alpha}R_{\alpha}^Af=0\quad,\alpha>0,f\in\mathcal
B^+(\overline{\Omega})$$ If $\mu\in S_{00}(\po)$, and if $f$ is
bounded function in $L^2(\Omega)$, then $||R_{\alpha}Â
f||<\infty$, and $U_A^{\alpha}R_{\alpha}^Af$ is a relative quasi
continuous version of the $\alpha-$potential
$U_{\alpha}(R_{\alpha}^Af.\mu)\in\mathcal F$ by Proposition 2.6.
Since
$||U_{\alpha}(R_{\alpha}^Af.\mu)||_{\infty}\leq||R_{\alpha}^Af||_{\infty}||U_{\alpha}\mu||_{\infty}<\infty$
and $\mu(\partial\Omega)<\infty$, we have that
$$R_{\alpha}^Af=R_{\alpha}f-U_A^{\alpha}R_{\alpha}^Af\in\mathcal
F^{\mu}$$ and that

\begin{equation*}
\begin{split}
\mathcal E_{\alpha}(R_{\alpha}^Af,u)& = \mathcal E_{\alpha}(R_{\alpha}f,u)-\mathcal E_{\alpha}(U_A^{\alpha}R_{\alpha}^Af,u)\\
                  &
                  =(f,u)-(R_{\alpha}^Af,u)_{\mu}\qquad,u\in\mathcal
                  F^{\mu}
\end{split}
\end{equation*}
$(4.2)$ follows.

For general positive measure $\mu$ charging no set of zero
relative capacity, we can take by virtue of Theorem 2.3 and Lemma
2.2 an increasing sequence $(F_n)$ of generalized nest of
$\partial \Omega$, and $\mu_n=1_{F_n}.\mu\in S_{00}(\po)$. Since
$\mu$ charges no set of zero relative capacity, $\mu_n(B)$
increases to $\mu(B)$ for any $B\in\mathcal B(\partial\Omega)$.

Let $A_n=1_{F_n}.A$. Then $A_n$ is a PCAF of $X_t$ with Revuz
measure $\mu_n$. Since $\mu_n\in S_{00}(\po)$ we have for $f\in
L^2(\Omega)$:
\begin{equation}
R_{\alpha}^{A_n} f\in\mathcal F^{\mu_n}\quad,\mathcal
E_{\alpha}^{\mu_n}(R_{\alpha}^{A_n},u)=(f,u)\quad,f\in
L^2(\Omega,m),u\in\mathcal F^{\mu_n}
\end{equation}
Clearly
$|R_{\alpha}^{A_n}f|\leq R_{\alpha}|f|<\infty$ r.q.e, and hence
$\lim_{n\rightarrow+\infty}R_{\alpha}^{A_n}f(x)=R_{\alpha}^Af(x)$
for r.q.e $x\in\overline{\Omega}$. For $n<m$, we get from $(4.3)$:
\begin{equation}
\mathcal
E_{\alpha}^{\mu_n}(R_{\alpha}^{A_n}f-R_{\alpha}^{A_m}f,R_{\alpha}^{A_n}f-R_{\alpha}^{A_m}f)\leq(f,R_{\alpha}^{A_n}f-R_{\alpha}^{A_m}f)
\end{equation}
which converges to zero as $n,m\rightarrow +\infty$. Therefore
$(R_{\alpha}^{A_n}f)_n$ is $\mathcal E_1-$convergent in $\mathcal
F$ and the limit function $R_{\alpha}^Af$ is in
$\widetilde{\mathcal F}$. On the other hand we also get from
$(4.3)$:

$
||R_{\alpha}^{A_n}f||_{L^2(\po,\mu)}\leq(f,R_{\alpha}^{A_n}f)_{L^2(\Omega)}\leq\frac{1}{\alpha}||f||_{L^2(\Omega)}$.
And by Fatou's Lemma:
$||R_{\alpha}^{A}f||_{L^2(\Omega)}\leq\frac{1}{\sqrt{\alpha}}||f||_{L^2(\Omega)}$,
getting $R_{\alpha}^Af\in\mathcal F^{\mu}$. Finally, observe the
estimate: $$
|(R_{\alpha}^{A_n}f,u)_{\mu_n}-(R_{\alpha}^{A}f,u)_{\mu}
|\leq||R_{\alpha}^{A_n}f-R_{\alpha}^{A}f||_{L^2(\partial\Omega,\mu_n)}||u||_{L^2(\partial\Omega,\mu)}+|(R_{\alpha}f,u)_{\mu-\mu_n}|$$holding
for $u\in L^2(\partial\Omega,\mu)$. The second term of the
right-hand side tends to zero as $n\rightarrow +\infty$. The first
term also tends to zero because we have from $(4.3)$: $\quad
||R_{\alpha}^{A_n}f-R_{\alpha}^{A_m}f||_{L^2(\partial\Omega,\mu_n)}\leq(f,R_{\alpha}^{A_n}f-R_{\alpha}^{A_m}f)$,
and it suffices to let first $m\rightarrow +\infty$ and then
$n\rightarrow +\infty$. By letting $n\rightarrow +\infty$ in
$(4.2)$ we arrive to desired equation $(4.1)$.

\end{proof}

The proof of the Theorem 4.2 is similar to the Theorem 6.1.1
\cite{FOT} which was formulated in the first time by S. Albeverio
and Z. M. Ma \cite{AM1} for general smooth measures in the context
of general $(X,m)$. In the case of $X=\overline{\Omega}$, and
working just with measures on $S_0(\partial\Omega)$ the proof
still the same, and works also for any smooth measure concentrated
on $\partial\Omega$. Consequently, the theorem still verified for
smooth measures ''nowhere Radon'' i.e. measures locally infinite
on $\partial\Omega$.

\begin{example} We give some particular examples of $\pmu$:

(1) If $\mu=0$, then $$\mathcal P_t^0 f(x)=E_{x}[f(X_t)] $$ the
semigroup generated by Laplacian with Neumann boundary conditions.

(2) If $\mu$ is locally infinite (nowhere Radon) on $\po$, then
$$\mathcal P_t^{\infty}f(x)=E_{x}[f(B_t)1_{\{t<\tau\}}]$$the semigroup
generated by the Laplacian with Dirichlet boundary conditions (see
Proposition 3.2.1 \cite{Wa}).

(3)  Let $\Omega$ be a bounded and enough smooth to insure the
existence of the surface measure $\sigma$, and $\mu=\beta.\sigma$,
with $\beta$ is a measurable bounded function on $\po$, then
$A_t^{\mu}=\int_0^t\beta(X_s)dL_s$, where $L_t$ is a boundary
local time. Consequently :$$\pmu
f(x)=E_x[f(X_t)exp(-\int_0^t\beta(X_s)dL_s)]$$is the semigroup
generated by the Laplacian with (classical) Robin boundary
conditions given by (1.1).
\end{example}
\bigskip

The setting of the problem from the stochastic point of view and the stochastic representation of the solution of the problem studied are important on themselves and are new. In fact  before there was always additional hypothesis on the domain or on the class of measures. Even if our approach is inspired by the works \cite{AM1}, \cite{AM2} and chapter 6 of \cite{FOT}, the link is not obvious and give arise to a new approach to the Robin boundary conditions. As a consequence, the proof of many propositions and properties become obvious and direct.The advantage of the stochastic approach is then, to give explicitly the representation of the semigroup and an easy access of it.

\begin{proposition}
$\mathcal P_t^{\mu}$  is sub-markovian i.e. $ \pmu\geq 0$
for all $t\geq 0$, and $$||\pmu
f||_{\infty}\leq||f||_{\infty}\quad (t\geq 0)$$
\end{proposition}

\begin{proof}
It is clear that if $f\in L^2(\Omega)_+$, then $\pmu f\geq 0$ for
all $t\geq 0$. In addition we have: $|\pmu f(x)|\leq
E_x[|f|(X_t)]$, and then $||\pmu
f||_{\infty}\leq||f||_{\infty}\quad (t\geq 0)$
\end{proof}

\begin{remark}
The analytic proof need the first and the second Beurling–Deny criterion (Proposition 3.10 \cite{AW1}) while  our proof is obvious and direct.
\end{remark}

Let $\Delta_{\mu}$ be the self-adjoint operator on $L^2(\Omega)$
generator of the semigroup $\mathcal P_t^{\mu}$, we write:
$$\mathcal P_t^{\mu}f(x)=e^{-t\Delta_{\mu}}f(x)$$ Following
\cite{Wa}, we know that $\Delta_{\mu}$ is a realization of the
Laplacian. Then we call $\Delta_{\mu}$ the Laplacian with General
Robin boundary conditions.

\begin{theorem}
Let $\mu\in S(\po)$, then the semigroup $\mathcal P_t^{\mu}$ is
sandwiched between the semigroup of Neumann Laplacian, and the
semigroup of Dirichlet Laplacian. That is :$$0\leq
e^{-t\Delta_D}\leq\mathcal P_t^{\mu}\leq e^{-t\Delta_N}$$for all
$t\geq 0$, in the sense of positive operators.
\end{theorem}

\begin{proof}

Let $f\in L^2(\Omega)_{+}$. Since $A_t^{\mu}\geq 0$
we get easily: $\mathcal P_t^{\mu}f(x)\leq
E_x[f(X_t)]$ for any $x\in\overline{\Omega}$. In the other hand we
have: $\mathcal P_t^{\mu}f(x)\geq
E_x[f(X_t)e^{-A_t^{\mu}}1_{\{t<\sigma_{\po}\}}]$, where
$\sigma_{\po}$ is the first hitting time of $\po$. Since the
relative quasi-support of $A_t^{\mu}$ and $N_t$ are in
$\partial\Omega$, then in $\{t<\sigma_{\po}\}$, $N_t$ and
$A_t^{\mu}$ vanishes. Consequently $X_t=B_t$ in
$\{t<\sigma_{\po}\}$ and $\mathcal P_t^{\mu}f(x)\geq
E_x[f(B_t)1_{\{t<\sigma_{\po}\}}]$. The theorem follows.
\end{proof}

\begin{remark}
The fact that the semigroup $\mathcal P_t^{\mu}$ is sandwiched between the Neumann semigroup and the Dirichlet one as proved in \cite{Wa} (Theorem 3.4.1) is not obvious and need a result characterizing domination of positive semigroups due to Ouhabaz, while our proof is simple and direct.
\end{remark}

\begin{proposition}Let $\mu,\nu\in S(\partial\Omega)$ such that $\nu\leq\mu$
(i.e $\nu(A)\leq\mu(A),\forall A\in\mathcal B(\partial\Omega)$),
then
$$0\leq e^{-t\Delta_D}\leq\mathcal P_t^{\mu}\leq\mathcal
P_t^{\nu}\leq e^{-t\Delta_N}$$for all $t\geq 0$, in the sense of
positive operators.
\end{proposition}

\begin{proof}It follows from the remark that if $\nu\leq\mu$, then
$A_t^{\nu}\leq A_t^{\mu}$, which means that $(A_t^{\mu})_{\mu}$ is
increasing, and then $(\mathcal P_t^{\mu})_{\mu}$ is decreasing.

\end{proof}

There exist a canonical Hunt process $X_t^A$ possessing the
transition function $\mathcal P_t^{\mu}$ which is directly
constructed from $X_t$ by killing the paths with rate $-dL_t$,
where $L_t=e^{-A_t}$.

To construct the process associated with $\mathcal P_t^{\mu}$, we
follow A.2 of \cite{FOT}, so we need a nonnegative random variable
$Z(\omega)$ on $(\Xi,\mathcal M, P_x)$ which is of an exponential
distribution with mean $1$, independent of $(X_t)_{t\geq 0}$ under
$P_x$ for every $x\in\ombar$ satisfying
$Z(\theta_s(\omega))=(Z(\omega)-s)\vee 0$. Introducing now a
Random time  $\xi^A$ defined by:  $$\xi^A=inf\{t\geq 0:A_t\geq
Z\}$$

We define the process $(X_t^A)_{t\geq 0}$ by:

$$X_t^A=\left\{
         \begin{array}{ll}
          X_t &  \text{if} \quad t<\xi^A; \\
           \Delta &
           \text{if} \quad t\geq \xi^A
         \end{array}
         \right.$$where $\Delta$ is a one-point compactification.

And, the admissible filtration of the process $(X_t^A)_{t\geq 0}$
is defined by:

$$\mathcal F_t^A=\{\Lambda\in\mathcal F_{\infty}:\Lambda\cap\{A_t<Z\}=\Lambda_t\cap\{A_t<Z\},\exists\Lambda_t\in\mathcal
F_t\}$$

Since $\{A_t<Z\}\cap\{A_t=\infty\}=\emptyset$, we may and shall
assume that $\Lambda_t\supset\{A_t=\infty\}$.

Now, we can write:

\begin{equation}
\begin{split}
E_x[f(X_t^A)]& = E_x[f(X_t):t<\xi^A]\\
            & = E_x[f(X_t):A_t<Z]\\
            & = E_x[f(X_t)e^{-A_t}]\\
            & = \pmu f(x)
\end{split}
\end{equation}

The Hunt process $(X_t^A)_{f\geq 0}$ is called the canonical
subprocess of $(X_t)_{t\geq 0}$ relative to the multiplicative
functional $L_t$. In fact, $(X_t^A)_{t\geq 0}$ is a Diffusion
process as $(\me^{\mu},\mf^{\mu})$ is local.

In the literature the Diffusion process $X_t^A$ is called
Partially reflected Brownian motion \cite{Gre}, in the sense that,
the paths of $X_t$ are reflected on the boundary since they will
be killed (absorbed) at the random time $\xi^A$ with rate $-d
L_t$.

\begin{theorem}Let $\mu,\mu_n\in S(\po)$ such that $\mu_n$ is monotone and
converges setwise to $\mu$ i.e.,$\mu_n(B)$ converges to $\mu(B)$
for any $B\in\mathcal B(\po)$, then $\Delta_{\mu_n}$ converges to
$\Delta_{\mu}$ in strongly resolvent sense.
\end{theorem}

\begin{proof}
We prove the theorem for $\mu_n$ increasing, the proof of the
decreasing case is similar. Let $A_n$ (respectively A) be the
additive functional associated to $\mu_n$( respectively $\mu$) by
the Revuz correspondence. Similarly to the second part of the
proof of Theorem 4.2, we have $\lim
R_{\alpha}^{A_n}f(x)=R_{\alpha}^Af(x)$ for r.q.e $x\in\ombar$.
Consequently $\lim_{n\rightarrow
+\infty}||R_{\alpha}^{A_n}f-R_{\alpha}^{A}f||_{L^2(\Omega)}=0$.
For $n<m$, we have $\mf^{\mu_m}\subset\mf^{\mu_n}$, and then
$$\quad\mathcal
E_{\alpha}^{\mu_n}(R_{\alpha}^{A_n}f-R_{\alpha}^{A_m}f,R_{\alpha}^{A_n}f-R_{\alpha}^{A_m}f)\leq(f,R_{\alpha}^{A_n}f-R_{\alpha}^{A_m}f)$$which
converges to zero as $n,m\rightarrow +\infty$. Therefore
$(R_{\alpha}^{A_n}f)_n$ is $\mathcal E_1-$convergent in $\mathcal
F$ and the limit function $R_{\alpha}^Af$ is in
$\widetilde{\mathcal F}$. The result follows.

\end{proof}

\begin{corollary}
Let $\mu\in S(\po)$ finite and let $k\in\mathbb{N}^*$. We defined
for $u,v\in\mathcal F^{\mu}$:$$\mathcal
E^{\mu_k}(u,v)=\int_{\Omega}\nabla u \nabla v dx
+\frac{1}{k}\int_{\po}\widetilde{u}\widetilde{v}d\mu$$ then
$\Delta_{\mu_k}\rightarrow\Delta_{N}$ in the strong resolvent
sense.

\end{corollary}

Intuitively speaking, when the measure $\mu$ is
infinity (locally infinite on the boundary), the semigroup $\mathcal P_t^{\mu}$ is the Dirichlet semigroup as said in the example 2, which mean that the boundary became "completely absorbing", and any
other additive functional in the boundary can not influence this
phenomena, which explain why $N_t$ doesn't appear yet in the
decomposition of $X_t$, which means that the reflecting phenomena
disappears, and so any path of $X_t$ is immediately killed when it
arrives to the boundary.

When $\mu$ is null on the boundary, then the semigroup $\mathcal P_t^{\mu}$ is the Neumann one, and then the boundary became ‘’completely reflecting’’, but for a general measure $\mu$ the paths are reflected many times before they will be absorbed at a random
time.


\par\bigskip



\begin{thebibliography}{99}
\addcontentsline{toc}{chapter}{Bibliographie}

\bibitem{AM1} Albeverio S. and Ma Z.-M.: Perturbation of Dirichlet forms-lower semiboundedness, closability,
and form cores. {\em J. Funct. Anal.} {\bf 99} (1991), 332-356.

\bibitem{AM2} Albeverio S.and Ma Z.-M.: Additive functionals, nowhere Radon and Kato class smooth
measures associated with Dirichlet forms. {\em Osaka J. Math},
{\bf 29} (1992) ), 247-265.

\bibitem{AW1}  Arendt W. and Warma M.: The Laplacian with Robin boundary conditions
on arbitrary domains. {\em Potential Anal.} {\bf 19} (2003),
341-363.

\bibitem{AW2} Arendt W. and Warma M.: Dirichlet and Neumann boundary conditions:
What is in between?.{\em J.evol.equ.} {\bf 3} (2003) 119 – 135.

\bibitem{BBC} Bass R. F. , Burdzy  K. and Chen Z-Q.: On the Robin problem in fractal
domains.{\em Proc. London Math. Soc.} Page 1 of 39 (2007).

\bibitem{BBC2} Bass R. F., Burdzy  K. and Chen Z.-Q.: Uniqueness for reflecting Brownian motion in lip domains.{\em Ann.
Inst. H. Poincare} {\bf 41} (2005) 197-235.

\bibitem{BH2} Bass  R.F. and Hsu E. P.: The semimartingale structure of reflecting Brownian
motion.{\em Proceedings of the American Mathematical Society},
{\bf 108}, (1990) 1007-1010.

\bibitem{BH} Bass R.F.
and Hsu E.P.: Some potential theory for reflecting brownian motion
in Holder and Lipschiz domains.{\em The Annals of Probability},
Vol 19, No 2, (1991), 486-508.

\bibitem{BM} Blanchard  Ph. and Ma Z.-M: New Results on the Schrodinger Semigroups with Potentials given
by Signed Smooth Measures. preprint,

\bibitem{CFW} Chen Z.Q., Fitzsimmons P.J. and Williams R.J.: Quasimartingales and strong Caccioppoli
set. {\em Potential Analysis}, {\bf 2} (1993), 281-315.

\bibitem{F} Fukushima M.: A construction of reflecting barrier Brownian motions for bounded
domains, {\em Osaka J. Math.}, {\bf 4} (1967), 183-215.

\bibitem{F2} Fukushima M.: Dirichlet forms, Caccioppoli sets and the Skorohod
equations, {\em Stochastic differential and difference equations}, (1997), 59-66, MR1636827.

\bibitem{FOT} Fukushima M., Oshima Y. and Takeda M.: Dirichlet Forms and Symmetric Markov Processes.
{\em Walter de Gruyter}, Berlin, (1994).


\bibitem{Gre}  GrebenkovD. S.: Partially Reflected Brownian Motion: A Stochastic Approach
to Transport Phenomena. arXiv:math.PR/0610080 v1 2 Oct (2006).

\bibitem{H3} Hsu E. P.: Reflecting Brownian motion, boundary local time and the Neumann
problem. Thesis, June 1984, Stanford University.

\bibitem{K} Kato T.: Perturbation theory for linear
operators.{\em  Classics In Mathematics}, Springer-Verlag Berlin
Heidelberg (1995) .

\bibitem{P} Papanicolaou V. G.: The probabilistic solution of the third boundary value problem for second order
elliptic equations. {\em Probab. Theory Related Fields}, {\bf 87}
(1990). 27-77.

\bibitem{R} Ramasubramanian S.: Reflecting brownian motion in
a lipschiz domain and a conditional gauge theorem.{\em Sankhya :
The Indian Journal of Statistics} Volume 63, Series A, Pt. 2,
(2001) pp. 178-193.

\bibitem{Re} Renming S.: Probabilistic approach to the third boundary value problem. Preprint, (1987).

\bibitem{SU} Sato K. and Ueno T.: Multi-dimensional diffusion and the Markov
process on the boundary.{\em J. Math. Kyoto Univ.} {\bf 4}-3
(1965) 529-605.

\bibitem{Si} Simon B.:t Schrodinger semigroups.{\em Bull. Amer.
Math. Soc} (N.S) {\bf 7} (1982) 447-526,

\bibitem{S} Stollmann P.: Smooth perturbation of Regular Dirichlet
Forms.{\em Proc. of the Amer. Math. Soc.}, Vol.116, No. 3.
(Nov.,1992), pp. 747-752.

\bibitem{SV} Stollman P. and Voigt J.: Perturbation of Dirichlet forms by
measures.{\em Potential Analysis}, {\bf 5} (1996), 109-138.

\bibitem{V1} Voigt J.: Absorption semigroups.{\em J. Operator
Theory}, {\bf 30} (1988), 117-131.

\bibitem{V2} Voigt J.: Absorption Semigroups, their Generators,
and Schrodinger Semigroups.{\em J. Functional Analysis,} {\bf 67},
(1986) 167-205 .

\bibitem{Wa} Warma M.: The Laplacian with General
Robin Boundary Conditions. Ph.D. thesis, University of Ulm,
(2002).




\end{thebibliography}
\end{document}